\documentclass[preprint,12pt,letterpaper]{elsarticle}

\usepackage[margin=1in]{geometry}
\usepackage[T1]{fontenc}
\usepackage[utf8]{inputenc}
\usepackage{lmodern}
\usepackage[english]{babel}
\usepackage{amsmath,amssymb,amsthm}
\usepackage{microtype}
\usepackage[hidelinks]{hyperref}
\usepackage{cleveref}

\theoremstyle{plain}
\newtheorem{thm}{Theorem}[section]
\newtheorem{lem}{Lemma}[section]
\newtheorem{cor}{Corollary}
\theoremstyle{definition}
\newtheorem{defn}{Definition}
\newtheorem{rem}[thm]{Remark}
\Crefname{thm}{Theorem}{Theorems}
\crefname{lem}{lemma}{lemmas}
\Crefname{lem}{Lemma}{Lemmas}

\journal{Discrete Mathematics}
\hypersetup{
  pdftitle={An Exact Vertex Localization of Zykov's Theorem},
  pdfauthor={Rajat Adak; L. Sunil Chandran},
  pdfsubject={Discrete Mathematics -- Note},
  pdfkeywords={Zykov's theorem, Turan's theorem, vertex-based localization, clique counts}
}

\begin{document}
\begin{frontmatter}
\title{An Exact Vertex Localization of Zykov's Theorem}
\author[iisc]{Rajat Adak\corref{cor1}}
\ead{rajatadak@iisc.ac.in}
\author[iisc]{L. Sunil Chandran}
\ead{sunil@iisc.ac.in}
\cortext[cor1]{Corresponding author}
\affiliation[iisc]{organization={Department of Computer Science and Automation,
Indian Institute of Science}, city={Bengaluru}, country={India}}
\begin{abstract}
We establish a vertex-based localization of the exact
form of Zykov's theorem. For a graph $G$ of order $n\geq1$,
let $N(G,K_t)$ denote its number of copies of $K_t$, and
let $c(v)$ be the order of a largest clique containing
$v$. Write $Z_t(n,s)=N(T_s(n),K_t)$, where $T_s(n)$ is
the balanced complete $s$-partite graph on $n$ vertices.
We prove that, for every integer $t\geq2$,
\[
N(G,K_t)\leq
\frac1n\sum_{v\in V(G)}Z_t(n,c(v)).
\]
Equality holds if and only if $G$ is $K_t$-free or
$G\cong T_s(n)$ for some $t\leq s\leq n$.
Moreover, the coefficients $Z_t(n,s)$ are pointwise
best possible among all bounds of this form.

Our result strengthens the previous vertex-localized
bound
\[
N(G,K_t)\leq n^{t-1}
\sum_{v\in V(G)}
\frac{\binom{c(v)}{t}}{c(v)^t}
\]
by accounting for the integer part sizes of Tur\'an
graphs. The comparison between the two upper bounds
is strict precisely when some vertex $v$ satisfies
$c(v)\geq t$ and $c(v)\nmid n$. As consequences, we
recover the exact Zykov and Tur\'an bounds without
divisibility assumptions.
\end{abstract}
\end{frontmatter}

\section{Introduction}

Extremal graph theory is concerned with maximizing or minimizing
graph parameters subject to given structural constraints, as well
as characterizing the graphs that attain the corresponding
extremal values. Let $m$ and $n$ denote the size and order of a
graph, respectively. One of the fundamental results in extremal
graph theory is the following classical theorem.

\begin{thm}\label{turan_main}
\emph{(Tur\'an \cite{Turan})}
Let $1\leq r\leq n$, and let $G$ be a $K_{r+1}$-free graph
on $n$ vertices. Then
\[
m\leq |E(T_r(n))|,
\]
where $T_r(n)$ is the Tur\'an graph on $n$ vertices with
$r$ parts. Equality holds if and only if $G\cong T_r(n)$.
\end{thm}

Let $N(G,K_t)$ denote the number of copies of the clique
$K_t$ in the graph $G$. Tur\'an's theorem was later
generalized by Zykov~\cite{Zykov1949-qr} in 1949, who
established an upper bound on $N(G,K_t)$ for
$2\leq t\leq r$ in any $K_{r+1}$-free graph $G$.

\begin{thm}\label{zykovthm}
\emph{(Zykov \cite{Zykov1949-qr})}
Let $2\leq t\leq r\leq n$, and let $G$ be an $n$-vertex
$K_{r+1}$-free graph. Then
\[
N(G,K_t)\leq N(T_r(n),K_t),
\]
with equality if and only if $G\cong T_r(n)$.
\end{thm}

For $2\leq t\leq r\leq n$, note that
$|E(T_r(n))|\leq \frac{n^2(r-1)}{2r}$ and
$N(T_r(n),K_t)\leq \binom{r}{t}(\frac{n}{r})^t$.
Equality in either inequality holds if and only if
$r\mid n$. These simpler expressions give the following
numerical relaxations of Tur\'an's and Zykov's theorems,
which are often used in the literature.

\begin{thm}\label{relaxbounds}
Let $2\leq t\leq r\leq n$, and let $G$ be an $n$-vertex
$K_{r+1}$-free graph with $m$ edges. Then
\[
m\leq\frac{n^2(r-1)}{2r}
\qquad\text{and}\qquad
N(G,K_t)\leq\binom{r}{t}\left(\frac nr\right)^t.
\]
Equality in either bound holds if and only if
$G\cong T_r(n)$ and $r\mid n$.
\end{thm}

For $2\leq t\leq r\leq n$, the numerical relaxation of
Zykov's bound in \Cref{relaxbounds} is exact when $r\mid n$
and strictly weaker otherwise. Indeed, if $n=ar+b$ with
$0\leq b<r$, the $r$ part sizes of $T_r(n)$ are $a+1$
($b$ times) and $a$ ($r-b$ times). Consequently, the
right side of the exact bound is
\begin{equation}
N(T_r(n),K_t)
=\sum_{j=0}^{t}
\binom{b}{j}\binom{r-b}{t-j}(a+1)^j a^{t-j},
\end{equation}
where an inadmissible binomial coefficient is interpreted
as zero. Thus, the exact bound retains information about
the integer part sizes that is lost in the numerical
relaxation. The aim of this paper is to retain this exact
dependence on the integer part sizes within the
vertex-based localization framework.

\subsection{Localization}

Recently, Brada\v{c}~\cite{bradac} and Malec and
Tompkins~\cite{DBLP:journals/ejc/MalecT23} introduced the
notion of \emph{localization} as a framework for
strengthening and extending classical extremal graph
bounds. Classical extremal results are often stated in
terms of global parameters such as clique number and
circumference.

The localization approach starts by defining a
\emph{localization function}. For instance, although the
clique number of a graph $G$ is a global parameter, it can
be localized by assigning to each edge $e\in E(G)$ a value
$w(e)$ equal to the size of the largest clique containing
$e$. Likewise, the circumference can be localized by
assigning to each edge the length of the longest cycle
containing it, with value $2$ for a cut-edge. Other global
parameters, including the length of a longest path,
maximum degree, and girth, can be localized in a similar
manner. Thus, localization replaces global parameters
with appropriate local weight functions defined on edges,
vertices, or other graph elements.

Brada\v{c}~\cite{bradac} and Malec and
Tompkins~\cite{DBLP:journals/ejc/MalecT23} gave an edge-based
localization of the numerical relaxation of Tur\'an's bound
as given in \Cref{relaxbounds}.

\begin{thm}\label{thm:GTuran}
\emph{(Brada\v{c} \cite{bradac} and Malec--Tompkins
\cite{DBLP:journals/ejc/MalecT23})}
Let $G$ be a graph on $n\geq1$ vertices. Then
\[
\sum_{e\in E(G)}\frac{w(e)}{w(e)-1}\leq\frac{n^2}{2}.
\]
Equality holds if and only if $G$ is a regular complete
multipartite graph with at least two nonempty parts.
\end{thm}

In \cite{adak2025vertex}, Adak and Chandran coined the term
\textit{vertex-based localization} for the framework to generalize extremal graph problems by assigning weights
to vertices. In \cite{adak2025turan}, the same authors
developed a vertex-based localization of the numerical
relaxation of Tur\'an's theorem.

There have been several developments in the localization
framework in recent years. In~\cite{zhao2025localized},
Zhao and Zhang used edge-based localization of the
circumference to generalize the Erd\H{o}s--Gallai theorem
for cycles~\cite{gallai1959maximal}, while the corresponding
theorem for paths was localized
in~\cite{DBLP:journals/ejc/MalecT23}.
Kirsch and Nir~\cite{Kirsch_2024} established localized
versions of several generalized Tur\'an-type results,
including a localization of the relaxed version of
Zykov's theorem. Arag\~{a}o and
Souza~\cite{aragao2024localised} further developed the
framework through localized versions of the graph
Maclaurin inequalities, using weight functions that
capture the local structure of the graph. Their results
also imply the localized bounds of Brada\v{c}, Malec and
Tompkins, and Kirsch and Nir. In \cite{adak2027vertex}, the authors obtained the vertex-based localized version of the generalized Erd\H{o}s-Gallai theorems.

Recently, a vertex-based localization of the numerical
relaxation of Zykov's theorem was presented
in \cite{adak2025generalizedzykovstheorem}. It also follows
from the weighted result of Arag\~{a}o and Souza
in \cite{aragao2024localised}.
Let $c:V(G)\rightarrow\mathbb{N}$ be defined by taking
$c(v)$ to be the order of a largest clique containing $v$.
Then, for every integer $t\geq2$,
\begin{equation}\label{weakzykov}
N(G,K_t)\leq n^{t-1}
\sum_{v\in V(G)}
\frac{1}{c(v)^t}\binom{c(v)}{t}.
\end{equation}

The coefficients in \eqref{weakzykov} still come from the
numerical relaxation of Zykov's bound. In particular, for
$G=T_s(n)$ with $t\leq s\leq n$ and $s\nmid n$, this upper
bound is strictly larger than $N(G,K_t)$. It is therefore
natural to seek a vertex-based localization that preserves
the exact extremal value for every graph order.

In this paper, we generalize \Cref{zykovthm} using the
vertex-based localization framework. Writing
$Z_t(n,s)=N(T_s(n),K_t)$, we prove that
\[
N(G,K_t)\leq
\frac1n\sum_{v\in V(G)}Z_t(n,c(v)).
\]
This bound strengthens \eqref{weakzykov} and is attained
by every $T_s(n)$ with $t\leq s\leq n$, without requiring
$s\mid n$. We characterize all equality cases and show
that the coefficients $Z_t(n,s)$ are best possible among
bounds of this form. As corollaries, we recover the exact
Zykov bound and obtain an exact vertex-based localization
of Tur\'an's theorem.

Our proof combines the transfer approach used
in \cite{adak2025generalizedzykovstheorem,aragao2024localised}
with an integer balancing argument. We also give examples
showing that the improvement over the classical bound
can be arbitrarily large, and that the same holds for
the improvement over \eqref{weakzykov} when the clique
size is allowed to vary.
\subsection{Preliminaries}
All graphs are finite and simple. For a positive integer $s$, write
$[s]=\{1,\ldots,s\}$.
\begin{defn}
For integers $s \in [n]$ and $t\ge 2$, set $Z_{t}(n,s):=N(T_{s}(n),K_{t})$. If $n=q_{s}s+\rho_{s}$, where $0\le\rho_{s}<s$, then
\begin{equation}\label{zt_defn}
Z_{t}(n,s) = \sum_{j=0}^{t}\binom{\rho_{s}}{j}\binom{s-\rho_{s}}{t-j}(q_{s}+1)^{j}q_{s}^{t-j} 
\end{equation}
In particular, $Z_{t}(n,s)=0$ when $s<t$, since there are no $t$-cliques in a Tur\'an graph with fewer than $t$ parts.
\end{defn}

\begin{defn}For nonnegative numbers $y_1,\ldots,y_s$, let
\begin{equation}\label{et_defn}
e_t(y_1,\ldots,y_s):=\sum_{1\le i_1<\cdots<i_t\le s}y_{i_1}\cdots y_{i_t}.\end{equation}
That is, $e_t(y_1,\ldots,y_s)$ is the sum of all products of $t$ coordinates with distinct indices from the tuple $(y_1,\ldots,y_s)$. We use the conventions $e_0(y_1,\dots,y_s)=1$, including
for the empty tuple, and $e_j(y_1,\dots,y_s)=0$ whenever $j>s$.\end{defn}

\begin{lem}\label{lem1}
Let $n\geq1$, $t\geq2$, and $1\leq s\leq n$ be integers.
Let $m_1,\dots,m_s$ be nonnegative integers with sum $n$.
Then
\begin{equation}
e_{t}(m_{1},\dots,m_{s})\le Z_{t}(n,s).
\end{equation}
If $2\le t\le s$ and all $m_{i}$ are positive, equality holds if and only if $|m_{i}-m_{j}|\le 1$ for every $i, j \in [s]$.
\end{lem}
    \begin{proof}
        Suppose there exists $i,j \in [s]$ such that $m_i \geq m_j +2$. Let $M= (m_1,\dots m_s)$ and $M'$ be the $s$-tuple obtained by replacing $m_i$ by $m_i -1$ and $m_j$ by $m_j +1$ in $M$. Clearly all the entries in $M'$ are nonnegative integers with sum $n$. Let $M_{-i,-j}$ denote the $s-2$-tuple of all the entries of $M$ except $m_i$ and $m_j$.
        Separating the terms according to whether they contain $i$ and $j$, we obtain
        \[
        e_t(M')-e_t(M)=(m_i-m_j-1)e_{t-2}(M_{-i,-j})\geq0.
        \]
        Each step decreases $\sum_{\ell=1}^{s}m_\ell^2$ by
        $2(m_i-m_j-1)>0$. Hence the process terminates at a tuple
        whose coordinates differ by at most one, proving the inequality.

         If all coordinates are positive and $2\leq t\leq s$, then $e_{t-2}(M_{-i,-j})>0$ whenever two entries differ by at least two (here $m_i$ and $m_j)$. Thus every step of transferring the values results in a vector with a strictly larger $e_t(\cdot)$ value, which proves the equality statement.
    \end{proof}

\begin{lem}\label{lem2}
For fixed $n, t$, the function $s\mapsto Z_{t}(n,s)$ is nondecreasing on $1\le s\le n$. Moreover,
\begin{equation}
Z_{t}(n,s+1)>Z_{t}(n,s) \quad (t\le s<n).
\end{equation}
\end{lem}
\begin{proof}
    By \Cref{lem1}, $Z_{t}(n,s)$ is the maximum of $e_{t}(m_{1},...,m_{s})$ over nonnegative integral $s$-tuples of sum $n$. Appending a zero gives a feasible $(s+1)$-tuple. Let $M =(m_1,\dots,m_s)$ be a balanced tuple, that is, $|m_i-m_j| \leq 1$ for every $i,j \in [s]$. Therefore, $Z_t(n,s) = e_t(M) = e_t(m_1,\dots,m_s,0)$. Thus, $Z_t(n,s+1) \geq e_t(m_1,\dots,m_s,0) = Z_t(n,s)$, proving monotonicity.
    
    For strictness, assume $t\leq s<n$. Every part of the balanced maximizing $s$-tuple is positive, and at least one part, say $m_i$, has size at least two. Append a zero and move one unit from such a part to the zero part. From \Cref{lem1}, the smoothing formula is strictly positive, that is,
    $e_t(m_1,\dots,m_i,\dots,m_s,0) < e_t(m_1,\dots,m_i-1,\dots,m_s,1).$
    
    Since, $Z_t(n,s) = e_t(m_1,\dots,m_i,\dots,m_s,0)$ and $Z_t(n,s+1) \geq e_t(m_1,\dots,m_i-1,\dots,m_s,1)$, we get the required strict inequality.
\end{proof}

\section{Our Result}
\begin{thm}\label{ourtheorem} Let $G$ be a graph of order $n\geq1$, let $t\geq2$ be an integer,
and define $c(v)$ as above. Then
\begin{equation} \label{eq:main}
N(G,K_{t})\le\frac{1}{n}\sum_{v\in V(G)}Z_{t}(n,c(v)).
\end{equation}
Equality holds if and only if one of the following holds:
(i) $G$ is $K_{t}$-free; or
(ii) $G\cong T_{s}(n)$ for some integer $s$ with $t\le s\le n$.
\end{thm}
\begin{rem}
Our bound strengthens \eqref{weakzykov}. Indeed, averaging the
part sizes over the reals, as in the proof of \Cref{lem1}, gives
\[Z_t(n,s)\leq\binom{s}{t}\left(\frac ns\right)^t.\]
For $t\leq s\leq n$, equality holds precisely when all parts
have equal size, equivalently, when $s\mid n$; for $s<t$,
both sides are zero. Consequently,
\[N(G,K_t)
\leq\frac1n\sum_{v\in V(G)}Z_t(n,c(v))
\leq n^{t-1}\sum_{v\in V(G)}
\frac{\binom{c(v)}{t}}{c(v)^t}.\]
The second inequality is strict if and only if some vertex
$v$ satisfies $c(v)\geq t$ and $c(v)\nmid n$.
In particular, for $G=T_s(n)$ with $t\leq s\leq n$ and
$s\nmid n$, our bound is attained, whereas the earlier
bound is strictly larger than $N(G,K_t)$.
\end{rem}
\begin{rem}The coefficients $Z_t(n,s)$ are best possible. Indeed, suppose
that a bound of the form
\[N(G,K_t)\leq\frac1n\sum_{v\in V(G)}f_t(n,c(v))\]
holds for every graph $G$ of order $n$. For any $s \in [n]$,
take $G=T_s(n)$. Since $c(v)=s$ for every vertex, we obtain
\[
Z_t(n,s)
=N(T_s(n),K_t)
\leq\frac1n\sum_{v\in V(T_s(n))}f_t(n,s)
=f_t(n,s).
\]
Thus, no coefficient $Z_t(n,s)$ can be decreased in a bound
of this form that holds for every graph of order $n$.\end{rem}
\subsection{Proof of \texorpdfstring{\Cref{ourtheorem}}{the main theorem}}
\begin{defn} Let $\Delta^{n-1}$ be the set defined as:
\[\Delta^{n-1} = \{x = (x_1,x_2,\dots,x_n) \in \mathbb{R}^n_{\geq 0} : nx_i \in \mathbb{Z}\ \forall\ i \in [n] \text{ and } \lVert x\rVert_1=1\}.\]
\end{defn}
\begin{defn}
    For a graph $G$ with $n$ vertices and a vector $x \in \Delta^{n-1}$, define $\Phi(G,x)=A(G,x)-B(G,x)$, where
    \[
    A(G,x)=\sum_{v\in V(G)}\frac{x_v}{n^t}Z_t(n,c(v)),
    \qquad
    B(G,x)=\sum_{Q\cong K_t}\prod_{v\in V(Q)}x_v.
    \]
    
Here each sum over $Q\cong K_t$ runs over the $t$-vertex cliques
of $G$, each counted once. We write
$\operatorname{supp}(x)=\{v\in V(G):x_v>0\}$.
\end{defn}
\begin{defn}
    We define a \textit{transfer function} on $\Delta^{n-1}$.
For $i, j \in V(G)$ and $x \in \Delta^{n-1}$, the transfer function transfers the amount $x_j$ from the $j$-th coordinate to the $i$-th coordinate of $x$. Formally,
\[
T^{ij}(x) = x + x_j(e_i - e_j).
\]
Here $e_i$ is the $i$-th standard unit vector. Note that it is easy
to verify that $T^{ij}(x)\in\Delta^{n-1}$.
\end{defn}

\begin{lem}\label{transfer}
    Since $\Delta^{n-1}$ is finite and nonempty, $\Phi(G,\cdot)$ attains a minimum on this set. Let $x$ minimize $\Phi(G,\cdot)$ on $\Delta^{n-1}$. If $i, j\in \operatorname{supp}(x)$ are nonadjacent, then each of the two vectors, $T^{ij}(x)$ and $T^{ji}(x)$, is also a minimizer.
\end{lem}
\begin{proof}
Let $y=T^{ij}(x)$ and $z=T^{ji}(x)$. For each vertex $v$, put
$$D_v(x):=\frac{Z_t(n,c(v))}{n^t}
-\sum_{\substack{Q\cong K_t\\v\in V(Q)}}
\prod_{w\in V(Q)\setminus\{v\}}x_w.$$

We first compare $\Phi(G,y)$ and $\Phi(G,x)$.
Only the coordinates of $i$ and $j$ change,
$y_i=x_i+x_j$ and $y_j=0$.
Since $i$ and $j$ are nonadjacent, no $t$-clique contains
both of them. Cliques containing neither vertex therefore
contribute the same amount to $B(G,y)$ and $B(G,x)$.
For a clique containing $i$, its contribution increases by
$x_j$ times the product of the other coordinates.
For a clique containing $j$, its contribution decreases by
$x_j$ times the product of the other coordinates.

Combining these changes with
\[
A(G,y)-A(G,x)
=\frac{x_j}{n^t}\bigl(Z_t(n,c(i))-Z_t(n,c(j))\bigr),
\]
we obtain
\[
\Phi(G,y)-\Phi(G,x)
=x_j\bigl(D_i(x)-D_j(x)\bigr).
\]
Similarly, transferring $x_i$ from $i$ to $j$ gives
\[
\Phi(G,z)-\Phi(G,x)
=x_i\bigl(D_j(x)-D_i(x)\bigr).
\]

Both differences are nonnegative because $x$ is a minimizer
and $y,z\in\Delta^{n-1}$. Since $x_i,x_j>0$, the first
identity implies $D_i(x)\geq D_j(x)$, while the second
implies $D_j(x)\geq D_i(x)$. Thus $D_i(x)=D_j(x)$,
and consequently
\[
\Phi(G,y)=\Phi(G,x)=\Phi(G,z).
\]
Therefore both transferred vectors are also minimizers.
\end{proof}
\subsubsection{Proof of Inequality}
Now we are all set to prove the inequality of \Cref{ourtheorem}.
\begin{proof}
    Choose a minimizer $x\in\Delta^{n-1}$ having minimum support. By \Cref{transfer}, two support vertices cannot be nonadjacent. Thus $G[\operatorname{supp}(x)]$ is a clique. Put $k=|\operatorname{supp}(x)|$. For $v\in \operatorname{supp}(x)$ we have $c(v)\geq k$. 
    For $v\in\operatorname{supp}(x)$, put $m_v=nx_v$.
    These are positive integers with sum $n$. Since the support
    is a clique, \Cref{lem1,lem2} give
    
    $$A(G,x)\geq\frac{Z_t(n,k)}{n^t}
    \geq\frac{e_t((m_v)_{v\in\operatorname{supp}(x)})}{n^t}
    =B(G,x) \implies A(G,x) - B(G,x) \geq 0.$$
    Thus $\Phi(G,x)\geq0$.
Evaluating at $u=(\frac{1}{n},\ldots,\frac{1}{n})\in\Delta^{n-1}$ gives
\[
0\leq\Phi(G,u)
=\frac1{n^t}\left(
\frac1n\sum_{v\in V(G)}Z_t(n,c(v))-N(G,K_t)
\right),
\]
which proves the inequality.
\end{proof}
\subsubsection{Proof of equality class}
\begin{proof}\textit{If part.} If $G$ is $K_{t}$-free, clearly $N(G,K_t) = 0$. Also, $c(v) < t$ for all $v \in V(G)$, therefore, $Z_t(n,c(v))=0$. Thus, both sides of the inequality are zero.

Now suppose that $G\cong T_s(n)$ for some $t\leq s\leq n$.
Every clique contains at most one vertex from each part.
Moreover, every vertex belongs to an $s$-clique, obtained by
choosing one vertex from each of the remaining parts.
Thus $c(v)=s$ for every vertex, and
\[\frac{1}{n}\sum_{v\in V(G)}Z_t(n,c(v))
=\frac1n\sum_{v\in V(G)}Z_t(n,s)
=Z_t(n,s)
=N(G,K_t).\]
Therefore equality holds in this case as well.

\medskip\noindent\textit{Only if part.} Now assume equality in the bound of \Cref{ourtheorem}. Therefore, the uniform vector $u$ must be a  minimizer of $\Phi(G,\cdot)$ in $\Delta^{n-1}$ and $\Phi(G,u) = 0$. 

Let $C_{1},...,C_{k}$ be the connected components of the complement $\overline{G}$ and let $p_{i}=|C_{i}|$. Choose an arbitrary root $r_{i}\in C_{i}$ and a rooted spanning tree of $\overline{G}[C_{i}]$. Starting at the leaves and proceeding toward the root, transfer the value of every non-root vertex to its parent, that is, if $i$ is the parent of $j$ in the spanning tree, then update the current vector, say $x$, to $T^{ij}(x)$. At each transfer, both the vertex and its parent still have positive coordinates, since neither has yet transferred its value to its own parent. Since the spanning tree is taken in $\overline{G}$, every tree edge is a non-edge of $G$. Thus, by repeated application of \Cref{transfer}, the resulting vector is a zero-valued minimizer $z$ with $z_{r_{i}}=\frac{p_{i}}{n}$ for $i \in [k]$, and $z_{v}=0$ elsewhere.

Distinct complement components have no complement edge between them; hence the roots form a $k$-clique in $G$, and $c(r_{i})\ge k$. Since $\Phi(G,z)=0$, we have $A(G,z)=B(G,z)$.
Therefore \Cref{lem1,lem2} give
\begin{equation}\label{chain}
\begin{aligned}
nZ_t(n,k)&\leq\sum_{i=1}^{k}p_iZ_t(n,c(r_i)) = n^{t+1}A(G,z)\\
&=n^{t+1}B(G,z) = n e_t(p_1,\ldots,p_k)\leq nZ_t(n,k).
\end{aligned}
\end{equation}
Thus equality holds throughout.

If $k<t$, then $Z_{t}(n,k)=0$. The equation implies $Z_{t}(n,c(r_{i}))=0$ for every $i \in [k]$, and therefore $c(r_{i})<t$. Each root was arbitrarily chosen, so $c(v)<t$ for every vertex. Hence $G$ is $K_{t}$-free.

If $k\geq t$, strict monotonicity in \Cref{lem2}, applied
to \eqref{chain}, yields $c(r_i)=k$ for every $i\in[k]$.
Since each root was chosen arbitrarily, $c(v)=k$ for every
$v\in V(G)$. The strict part of \Cref{lem1} yields
$|p_i-p_j|\leq1$ for all $i,j\in[k]$. There can be no edge
$uv$ of $G$ inside a component $C_j$, since $u,v$, together
with one vertex from each other component, would form a
$(k+1)$-clique. Thus, every $C_{i}$ is independent in $G$, while any two distinct $C_{i}$, $C_{j}$ are completely joined for $i,j \in [k]$. Therefore $G$ is complete $k$-partite, with parts differing by at most one. Hence $G\cong T_{k}(n)$.\end{proof}

\section{Retrieving Zykov's theorem}
\begin{cor}
Let $2\le t\le r\le n$, let $n=ar+b$ with $0\le b<r$ and let $G$ be an $n$-vertex $K_{r+1}$-free graph. Then
\begin{equation*}
N(G,K_{t})\le\sum_{j=0}^{t}\binom{b}{j}\binom{r-b}{t-j}(a+1)^{j}a^{t-j} = Z_t(n,r)=N(T_r(n),K_t)
\end{equation*}
Equality holds if and only if $G\cong T_{r}(n)$.
\end{cor}
\begin{proof}
Since $G$ is $K_{r+1}$-free, $c(v)\le r$ for every vertex. From \Cref{ourtheorem} and \Cref{lem2} we get
\[
N(G,K_{t})\le\frac{1}{n}\sum_{v\in V(G)}Z_{t}(n,c(v))\leq \frac{1}{n}\sum_{v\in V(G)}Z_{t}(n,r) = Z_{t}(n,r) = N(T_r(n),K_t).
\]
 If $G\cong T_{r}(n)$ equality is immediate. Conversely, if both ends are equal, equality holds in the localized bound. Since the clique count is positive, $G\cong T_{s}(n)$ for some $s\ge t$. Strict monotonicity forces $s=r$.
\end{proof}

\section{Comparison with the classical and previous localized bounds}\label{examples}

First, fix integers $2\leq t\leq r$. For each integer
$k\geq2$, let $G_k$ be obtained from $K_r$ by adding
$r(k-1)$ isolated vertices. Then $G_k$ has $n=kr$
vertices and clique number $r$. The classical Zykov
bound gives
\[
N(G_k,K_t)\leq Z_t(kr,r)=\binom{r}{t}k^t.
\]
The $r$ clique vertices have $c(v)=r$, while each isolated
vertex has $c(v)=1$ and contributes zero. Therefore,
\Cref{ourtheorem} gives
\[
N(G_k,K_t)
\leq\frac{1}{kr}\sum_{v\in V(G_k)}Z_t(kr,c(v))
=\frac{r}{kr}Z_t(kr,r)
=\binom{r}{t}k^{t-1}.
\]
The previous localized bound \eqref{weakzykov} also gives
$\binom{r}{t}k^{t-1}$, since $r\mid n$.
Thus, the exact localized bound improves the classical
bound by a factor of $k$, which can be arbitrarily large
even with $t$ and $r$ fixed.

To compare the two localized bounds, we now allow the
clique size to vary. For each integer $k\geq1$, let
$H_k=T_{2k}(3k)$ and take $t=2k$. The graph $H_k$ has
$k$ parts of size two and $k$ parts of size one.
Consequently, $c(v)=2k$ for every vertex, and
$N(H_k,K_{2k})=2^k$. Our exact localization gives
\[
N(H_k,K_{2k})
\leq\frac{1}{3k}\sum_{v\in V(H_k)}Z_{2k}(3k,c(v))
=Z_{2k}(3k,2k)
=2^k.
\]
Here the exact localized bound coincides with the classical
bound and is attained. In contrast, \eqref{weakzykov} gives
\[
N(H_k,K_{2k})
\leq(3k)^{2k-1}
\sum_{v\in V(H_k)}
\frac{\binom{c(v)}{2k}}{c(v)^{2k}}
=\left(\frac{9}{4}\right)^k.
\]
The previous localized upper bound is therefore
$(9/8)^k$ times the exact localized upper bound.
Since this ratio tends to infinity, the multiplicative
improvement over the previous localization is also
unbounded when the clique size is allowed to vary.
\section{An exact localization of Tur\'{a}n's theorem}

For $1\leq s\leq n$, write $n=q_s s+\rho_s$, where
$0\leq\rho_s<s$. Since $T_s(n)$ has $\rho_s$ parts of size
$q_s+1$ and $s-\rho_s$ parts of size $q_s$, we have
\begin{align*}
Z_2(n,s)
&=\frac12\left(n^2-\rho_s(q_s+1)^2-(s-\rho_s)q_s^2\right)\\
&=\frac{n^2(s-1)}{2s}-\frac{\rho_s(s-\rho_s)}{2s}.
\end{align*}
For each vertex $v$, let $b(v)$ be the remainder when $n$ is
divided by $c(v)$. Taking $t=2$ in \Cref{ourtheorem} gives
\begin{equation*}
\begin{aligned}
m&\leq\frac1n\sum_{v\in V(G)}Z_2(n,c(v))\\
&=\frac n2\sum_{v\in V(G)}\frac{c(v)-1}{c(v)}
-\frac1{2n}\sum_{v\in V(G)}
\frac{b(v)(c(v)-b(v))}{c(v)}.
\end{aligned}
\end{equation*}
By \Cref{ourtheorem}, equality holds if and only if
$G\cong T_s(n)$ for some $1\leq s\leq n$, including the edgeless
case $s=1$. The correction term vanishes precisely when
$c(v)\mid n$ for every vertex.

If $G$ is $K_{r+1}$-free and $n=ar+b$, where $0\leq b<r\leq n$,
then $c(v)\leq r$ for every vertex. Hence \Cref{lem2} recovers
the exact Tur\'{a}n bound
\[m\leq Z_2(n,r) = |E(T_r(n))|
=\frac{n^2(r-1)}{2r}-\frac{b(r-b)}{2r},\]
with equality if and only if $G\cong T_r(n)$.
\section{Conclusion}
We obtained an exact vertex localization of Zykov's theorem
that strengthens the earlier numerical localization by accounting
for the remainder in the part sizes. We characterized all equality
cases and recovered the exact Zykov and Tur\'an bounds for
arbitrary graph orders.
\section*{Declaration of generative AI and AI-assisted technologies}
ChatGPT was used solely to assist with correcting English grammar, rephrasing sentences, and improving the overall organization of the manuscript. The authors have reviewed and revised the resulting content and take full responsibility for the manuscript.
\section*{Acknowledgments}
Rajat Adak acknowledges support from the TCS Research Scholar Program (Cycle 20).
\bibliographystyle{plain}
\bibliography{references}
\end{document}